\documentclass{article}

 \usepackage{amsmath,amsthm,amssymb,amsfonts}
 \usepackage{indentfirst}

\usepackage{setspace}
\newcommand{\R}{\mathbb{R}}
\newcommand{\N}{\mathbb{N}}

\newcommand{\dive}{\mathrm{div}}
\newcommand{\F}{F}

\newenvironment{eremer}
               {\vspace{2mm}\noindent{\textbf{Acknowledgment}\nobreak}}
               {\vspace{2mm}}

\newtheorem{lem}{Lemma}
\newtheorem{prop}{Proposition}
\newtheorem{thm}{Theorem}
\newtheorem{cor}{Corollary}

\theoremstyle{remark}
\newtheorem{rem}{Remark}

\theoremstyle{remark}
\newtheorem{ex}{Example}
\theoremstyle{remark}

\newcommand{\PAR}[1]{\ensuremath{{\left(#1\right)}}} 

\begin{document}

\begin{center}

{\Large \bf Logarithmic Sobolev inequalities: regularizing effect of L\'evy operators and 
asymptotic convergence in the L\'evy-Fokker-Planck equation}
\bigskip

\textsc{Ivan Gentil}\footnotemark[1]
and  
\textsc{Cyril Imbert}\footnotemark[1]

\footnotetext[1]{
Centre de Recherche en Math\'ematiques de la D\'ecision 
(CEREMADE, UMR CNRS 7534), Universit\'e Paris-Dauphine, 
Place du Mar\'echal De Lattre De Tassigny
75775 Paris cedex 16, France, {\tt gentil,imbert@ceremade.dauphine.fr}}
{\today}

\end{center}

\begin{quote} \footnotesize
\noindent \textsc{Abstract.} 
In this paper we study some applications of the L\'evy logarithmic Sobolev inequality to the study of the regularity 
of the solution of the fractal heat equation, \textit{i. e.} the heat equation where the Laplacian is
replaced with the fractional Laplacian. It is also used to the study of  the asymptotic behaviour of the 
L\'evy-Ornstein-Uhlenbeck process. 
\end{quote}
\vspace{5mm}

\medskip


\section{Introduction}

On one hand, regularity results  for the heat equation in $\R^d$, such as ultracontractivity, can be obtained by using a Euclidean logarithmic Sobolev inequality. 
On the other hand the asymptotic behaviour of the Ornstein-Uhlenbeck semigroup, precisely the optimal exponential decay to the equilibrium, 
is proved by using  either Poincar\'e or logarithmic Sobolev inequalities. See \cite{logsob} for a review of the subject. 

The heat equation or the Ornstein-Uhlenbeck semigroup are associated with the Laplacian, the infinitesimal generator of the Brownian motion.
The Brownian motion makes part of a large class of stochastic processes called L\'evy processes.   
In this note we would like to describe how the properties we just mentioned (ultracontractivity and exponential decay)
are sometimes  true if we replace the Laplacian with the infinitesimal operator of a general L\'evy process. These 
generators are integrodifferential and are referred to as L\'evy operators.   

In the next section we  give a short introduction to L\'evy processes and L\'evy operators. 
Two important inequalities are also given:  the  Euclidean logarithmic inequality  in the case of the $\alpha$-stable process;
and a modified logarithmic Sobolev inequality for infinitely divisible probability measures; the latter inequality generalizes 
the logarithmic Sobolev inequality given by L. Gross in \cite{gross}.

In Section~\ref{sec-heat}, we prove that the heat equation associated with a $\alpha$-stable process satisfies the property of ultracontractivity. 

In Section~\ref{sec-fp} we consider the Ornstein-Uhlenbeck semigroup or equivalently the Fokker-Planck semi-group associated with a general 
L\'evy operator. We will see that, under proper assumptions on the operator, those semi-groups converge to the unique steady state. Results of  
Section~\ref{sec-fp} are presented in full details in~\cite{gi} and are extensions of the paper of P. Biler and G. Karch \cite{bk03}.

\section{Preliminaries}
\label{sec:prelim}

\subsection{L\'evy operators}
\label{subsec:levy}

Let us  recall basic definitions about L\'evy operators and introduce notations. 
See for example \cite{applebaum} for further details. 
\medskip

\noindent \textbf{Characteristic exponents and L\'evy measures.} Let $d\in \N^*$. 
 Because of the definition of a L\'evy process in $\R^d$ (a process with stationary and independent increment),
the law $\mu_t$ of such a process  $(X_t)_{t\geq0}$ at time $t>0$ is \emph{infinitely divisible}, 
\textit{i.e.} it can be written for all $n \ge 1$ under the form 
$$
\mu_t = \underset{\mbox{$n$ times}}{\underbrace{\mu_n \star \dots \star \mu_n}}
$$ 
for some probability $\mu_n$ (depending on $n$).  Using this property, it can be shown that the  characteristic function 
$\phi_{X_t} (\xi):=E(\exp(i\xi\cdot X_t))$ (\textit{i.e.} its Fourier transform)
of the law of $X_t$ can be written under the form $\exp ( t\psi(\xi))$ for
a function $\psi$ called the \emph{characteristic exponent}. The L\'evy-Khinchine formula states that $\psi$ can be described with exactly three parameters $(\sigma, b, \nu)$ where $\sigma$ is a nonnegative symmetric $d\times d$ matrix, $b\in\R^d$ and $\nu$ is a nonnegative singular measure on $\R^d$ that satisfies 
\begin{equation}
\label{eq-def1}
\nu (\{ 0\})=0 \quad \mbox{ and } \int \min (1, |z|^2) \nu (dz) <+\infty. 
\end{equation} 
Then $\psi$ can be written under the form
\begin{equation}
\label{def-psi}
\psi(\xi)= - \sigma \xi \cdot \xi + i b\cdot \xi  + a(\xi)
\end{equation}
where $a$ is given by 
$$
a(\xi)=\int \PAR{e^{iz\cdot\xi}-1-i(z\cdot\xi) h(z)}\nu(dz), 
$$
with  $h (z) = 1/({1+|z|^2})$. 

 The matrix $\sigma$ characterizes the diffusion (or Gaussian) part of the operator  (with eventually $\sigma=0$), while $b$
 characterizes the drift part and $\nu$ is called  a L\'evy measure; it characterizes
the pure jump part. The support of the measure $\nu$ represents the 
possible jumps of the process. 

A L\'evy operator $\mathcal{I}$ is the infinitesimal 
generator associated with  the L\'evy process and the L\'evy-Khinchine formula implies 
that it has the following form
\begin{equation}
\label{operator}
\mathcal{I} [u] (x) =\mathrm{div}\, (\sigma \nabla u)(x) +   b\cdot \nabla u(x) +
\int_{\R^d} (u(x+z) - u(x) - \nabla u (x) \cdot z  h  (z) ) \nu(dz)
\end{equation}

\medskip

\noindent \textbf{The pseudo-differential point of view.} It is convenient to introduce 
 the operator  $\mathcal I_g$ associated with the Gaussian part
$$\mathcal I_g(u)=\mathrm{div}\, (\sigma \nabla u) +   b\cdot \nabla u$$
 and the operator  $\mathcal I_a$ associated with the pure jump part
$$\mathcal I_a (u) = \int_{\R^d} ( u(x+z) -u (x) -\nabla u (x) \cdot z h(z)) \nu (dz).$$ 
The operator $\mathcal I_a$ can be seen as a pseudo-differential operator of 
symbol $a$
 $$
\quad{\mathcal I_a }(u) =  {\mathcal F} ( a \times {\mathcal F}^{-1} u )
 $$   
 where $\mathcal F$ stands for the Fourier transform (see Theorem~3.3.3 p.139 of \cite{applebaum}). Here, we choose the probabilistic
convention in defining,  for all function $w\in L^1(\R^d)$,
 \begin{equation}
 \label{eq-fourier}
 \forall \xi\in\R^d,\quad \hat{w}(\xi)={\mathcal F}(w)(\xi)=\int e^{i x\cdot \xi} w(x)dx. 
 \end{equation}

Moreover, using the Fourier interpretation of the L\'evy operator one gets the following integration by parts formula :
if $\mathcal{I}$ is a L\'evy operator with parameters $(b,\sigma,{\nu})$ then for any smooth functions $u$, $v$ one gets 
\begin{equation}
\label{eq-ipp}
\int v\mathcal{I}[u]dx=\int u \check{\mathcal{I}} [v]dx,
\end{equation}
where $\check{\mathcal{I}}$ is the L\'evy operator whose parameters are $(-b,\sigma,\check{\nu})$
 with $\check{\nu} (dz) = \nu (-dz)$.
 
 \medskip

\noindent \textbf{Multi-fractal and $\alpha$-stable L\'evy operators.}
L\'evy operators whose characteristic exponent is positively homogeneous of index $\alpha \in (0,2]$ 
are called \emph{$\alpha$-stable}. 
The fractional Laplacian corresponds to a particular $\alpha$-stable L\'evy process with characteristic exponent 
 $\psi (\xi) = |\xi|^\alpha$, where $|\cdot|$ is the Euclidean norm in $\R^d$.  In the case $\alpha \in (0,2)$, one gets
$$
b=0, \quad \sigma =0 \quad \mbox{  and } \quad \nu (dz) = \frac{dz}{|z|^{d+\alpha}}.
$$ 
Hence, 
it is a pure jump process, \textit{i.e.} it has neither a drift part nor a diffusion one. 

L\'evy operators whose characteristic exponent can be written  as $\psi(\xi)=\sum_{i=1}^n\psi_i(\xi)$ where $\psi_i$ is $\alpha_i$-homogenous  with $\alpha_i\in(0,2]$, are often referred to as multi-fractal L\'evy operators. 

\medskip

The $\alpha$-stable operators play a central role 
in this paper. Let $-g_\alpha [\cdot]$ denote the L\'evy operator associated to the $\alpha$-stable L\'evy process
whose characteristic exponent is 
 $\psi (\xi) = |\xi|^\alpha$, $\alpha\in(0,2]$.  In the limit case, when $\alpha=2$, $g_2 [\cdot]$ is the Laplacian operator on $\R^d$.

\subsection{A Euclidean logarithmic Sobolev inequatlity}
\label{subsec-eucli}

The  logarithmic Sobolev inequality for the Lebesgue measure
is a useful functional inequality in the study of the fractional heat equation: $\partial_t u+g_{\alpha}[u]=0$. 
Such an inequality has been established  by A. Cotsiolis and N. K. Tavoularis. 
\begin{thm}[\cite{ct}]
Let $\alpha\in(0,2]$ then for any smooth function $f$ on $\R^d$ such that $\int f^2 dx=1$, 
the following optimal Euclidean logarithmic Sobolev inequality holds true
\begin{equation}
\label{eq-lse}
{\rm{Ent}}_{dx}(f^2):=\int f^2\log f^2 dx\leq \frac{n}{\alpha}
\log\left(\frac{\alpha C^{\alpha/n}}{n\pi^{\alpha/2} e^{1-\alpha}}  \int(g_{\alpha/2}[f])^2 dx \right),
\end{equation}
where $C= \frac{2\Gamma(n/\alpha)}{\alpha\Gamma(n/2)}$.
\end{thm}

This inequality is a generalization of the classical Euclidean logarithmic Sobolev inequality given by F.B Weissler in \cite{wei80}. 

\subsection{A modified logarithmic Sobolev inequality}

In the sequel, we will need another functional inequality
 proved by C. An\'e and M. Ledoux \cite{ane-ledoux} 
in the particular case of the Poisson measure and then generalized 
by L. Wu \cite{wu00} and D. Chafa\"i \cite{chafai04} for all infinite measurable laws.
In order to state the most general result, we first introduce $\Phi$-entropies. 

Let  $\Phi: \R^+ \mapsto \R$ a smooth convex function and define the $\Phi$-entropy: for any nonnegative function $f$,
$$
\mathrm{Ent}^\Phi_{\mu}  \left( f \right)
: = \int \Phi \left( f \right) d\mu-\Phi\PAR{\int f d\mu}
$$
where $\mu$ is a probability measure. When $\Phi (x)=x\log x$ we recover the classical entropy introduced in~\eqref{eq-lse}.

For a convex function $\Phi$ we note by $D_\Phi$ the so-called {\it Bergman distance} defined by : 
\begin{equation}
\label{eq-breg}
\forall (a,b)\in\R^+,\quad D_\Phi (a,b) := \Phi (a) - \Phi (b) - \Phi'(b) (a- b) \ge 0.
\end{equation}

\begin{thm}[\cite{ane-ledoux,wu00,chafai04}]
\label{prop-ls}
Assume that $\Phi$ satisfies the following properties:
\begin{equation} \label{cond:convex}
\left\{\begin{array}{l}
(a,b) \mapsto D_\Phi (a+b,b) \\
(r,y) \mapsto \Phi''(r)  y \cdot\sigma y 
\end{array} \right. 
\quad \mbox{ are convex on } \quad 
 \left\{\begin{array}{l}
\{ a+b \ge 0, b \ge 0\} \\
\R^+ \times \R^{2d}
\end{array} \right. 
\end{equation}
where $D_\Phi$ has been defined by~\eqref{eq-breg}

Consider an infinitely divisible law $\mu$ on $\R^d$. 
Then for all smooth positive functions $v$,  
\begin{equation}
\label{eq-ls}
\mathrm{Ent}^\Phi_{\mu}  \PAR{v} \le 
\int\Phi''(v) \nabla v \cdot\sigma  \nabla v \mu (dx)
+ \int \int D_\Phi (v(x), v(x+z)) \nu_\mu (dz) \mu (dx)
\end{equation}
where $\nu_\mu$ and $\sigma$ denote respectively the L\'evy measure and the diffusion matrix 
associated with $\mu$. 
\end{thm}
Remark that the drift of the law plays no role in this functional inequality.  Inequality~\eqref{eq-ls} 
is proved in \cite{wu00} for $\Phi(x)=x^2$ or $\Phi(x)=x\log x$ and  in this general form in \cite{chafai04}. 
\medskip

An important special case is the following one: $\Phi (r) = r^2/2$. 
A simple computation shows that the Bregman distance $D_\Phi$ in this case is 
$D_\Phi (a,b) = (a-b)^2/2$ so that $\mathrm{Ent}^\Phi_{\mu}$ 
reduces to the variance (up to a constant). See also the appendix of \cite{gi} for a proof of
this inequality. 
 

\section{Regularity of the heat equation driven by a L\'evy process. }
\label{sec-heat}
In this section, we study regularity properties  of solutions of the fractional heat equation
\begin{equation}\label{eq:frac-heat}
\partial_t u+g_{\alpha}[u]=0.
\end{equation}
In particular, we are interested in the ultracontractivity of this equation. 
\begin{thm}
 \label{main}
Let $\alpha  \in (0,2)$ and   $(P_t)_{t\geq0}$ denote  the semigroup associated with the equation~\eqref{eq:frac-heat}. 
Consider a smooth initial datum $f$.  
Then  for all $t>0$ and $q \geq p \geq 2$
\begin{equation} \label{hc-diff1}
\| P_t f\|_q \le \| f \|_p \left( \frac{{\mathcal A}n (q - p)}{2
  \alpha t}\right)^{\frac{n (q -p)}{\alpha pq}}
  \frac{p^{n/(q\alpha)}}{q^{n/(p\alpha )}}
\end{equation}
where $\| \cdot \|_p$ denotes the $L^p(dx)$ norm and 
\begin{equation}\label{constanteA}
\mathcal A=\frac{\alpha \left(\frac{2\Gamma(n/\alpha)}{\alpha\Gamma(n/2)}\right)^{\alpha/n}}{n\pi^{\alpha/2}e^{1-\alpha}}.
\end{equation}
\end{thm}


These results can be found in the case of the Laplacian in \cite{bakry}. We would like to mention  that in the classical case, 
Inequality~\eqref{hc-diff1} for all $q\geq p\geq 2$ and the Euclidean logarithmic Sobolev inequality are equivalent which is not 
clear in our case.  
\medskip

 Letting $q \to +\infty$ and choosing $p=2$ in Theorem~\ref{main} yields:
\begin{cor}[Ultracontractivity]
The semigroup  $(P_t)_{t\geq0}$ is  \emph{ultracontractive}, i.e. it satisfies for all 
smooth function $f$
\begin{equation*}
\label{hc-diff}
 \| P_t f\|_\infty \le \|f \|_2 \left(\frac{{\mathcal A}n}{2\alpha
  t}\right)^{n/(2\alpha)}
\end{equation*}
where $\mathcal{A}$ is given by \eqref{constanteA}. 
\end{cor}
We next recall a useful inequality satisfied by L\'evy operator. Such an inequality
is sometimes called Kato inequality. See for instance the proof given in~\cite{di04}. 
\begin{lem}
\label{convexe}
Let $\phi: \R \to \R$ be convex and $u \in C^2_b(\R^N)$. Then,
if $\phi$ is differentiable at $u(x)$, we have:
$$ g_\alpha[\phi (u)] (x) \le \phi'(u(x)) g_\alpha [u](x).$$
\end{lem}
We will also use the simple fact that $\int u g_\alpha [v]dx = \int
|\xi|^\alpha \hat{u} \hat{v}dx$. In particular, for all smooth function $u$ on $\R^d$,
$$
  \int u g_\alpha [u]dx = \int(g_{\alpha/2} [u])^2dx.
$$ 
\begin{proof}[Proof of Theorem~\ref{main}]
Let $u(t,x)$ denote $P_t f (x)$ and consider an increasing function $\varphi:\R \to \R$
such that $\varphi(0)=\alpha$. Define a function $F(t)= \|u(t)\|_{\varphi(t)}$ and
let us study its derivative. A computation gives
\begin{multline*}
\frac{\varphi^2}{\varphi'} F^{\varphi-1} F' = \mathrm{Ent}_{dx} (|u|^\varphi) +
\frac{\varphi^2}{\varphi'}\int |u|^{\varphi-1}\partial_t u \, dx \\=  \mathrm{Ent}_{dx} (|u|^\varphi) -
\frac{\varphi^2}{\varphi'}\int |u|^{\varphi-1} g_\alpha [u] dx.
\end{multline*}
Assume that $\varphi\geq2$. In this case, one can apply Lemma~\ref{convexe} with $\phi(\cdot)=|\cdot|^{\varphi/2}$ and get
$$
 -\varphi |u|^{\varphi-1} g_\alpha [u] \le - 2 |u|^{\varphi/2} g_\alpha [|u|^{\varphi/2}];
 $$
integrating over $\R^n$ implies 
$$
-\varphi \int |u|^{\varphi-1} g_\alpha [u] dx \leq
-2 \int |u|^{\varphi/2} g_\alpha [|u|^{\varphi/2}] dx = -2 \int \left(g_{\alpha/2} [|u|^{\varphi/2}]\right)^2dx
$$ 
so that
$$
\frac{\varphi^2}{\varphi'} F^{\varphi-1} F' \leq \mathrm{Ent}_{dx} ((|u|^{\varphi/2})^2) -
\frac{2\varphi}{\varphi'}   \int \left(g_{\alpha/2} [|u|^{\varphi/2}]\right)^2dx.$$
Apply now \eqref{eq-lse} with $|u|^{\varphi/2} / \sqrt{\int |u|^\varphi dx}$ and get: 
\begin{multline*}
\frac{\varphi^2}{\varphi'} F^{\varphi-1} F' \leq
 \frac{n}{\alpha} \int |u|^\varphi dx \log\PAR{ {\mathcal A} 
\frac{ \int \left(g_{\alpha/2} [|u|^{\varphi/2}]\right)^2dx}{\int |u|^\varphi dx}} \\ 
-\frac{2\varphi}{\varphi'}  \int \left(g_{\alpha/2} [|u|^{\varphi/2}]\right)^2dx.
\end{multline*}
Use now the concavity of $\log$; for any $x\in \R$, we have
$$
\frac{\varphi^2}{\varphi'} F^{\varphi-1} F' \leq \left( \frac{n}{\alpha x } -
\frac{2\varphi}{\varphi'}\right) \int \left(g_{\alpha/2} [|u|^{\varphi/2}]\right)^2dx + \frac{n}{\alpha} \log
( {\mathcal A}x/e) \int |u|^\varphi dx .
$$
Choose next $x$ such that $2\varphi/\varphi' = n / (\alpha x) $. We now obtain that $F$ satisfies
$$ 
\frac{F'}{F} \le \frac{n\varphi'}{\alpha \varphi^2}  \log \left(
\frac{ {\mathcal A}n\varphi'}{2e\alpha \varphi} \right)
$$
so that for all $t>0$,
$$ 
\| P_t f\|_{\varphi(t)} = \| u(t)\|_{\varphi(t)} =  F(t) \le F(0) \exp \left\{ \int_0^t \frac{n \varphi'(s)}{\alpha \varphi^2(s)} 
\log  \left( \frac{ {\mathcal A}n\varphi'(s)}{2e\alpha \varphi(s)} \right) ds \right\}.
$$
We now minimize the right hand side of the previous inequality 
w.r.t. functions $\varphi$ such that $\varphi(0)=p$ and $\varphi(t)=q$. 
Associated Euler's equation reads: $2 \varphi'^2 =
\varphi'' \varphi$  so that we choose $\varphi(s) = \frac{t pq}{(p-q)s + q t}$ and one can check that with such
a choice of $\varphi$, Inequality \eqref{hc-diff1} is proved.

\end{proof}

\section{Asymptotic behaviour of a L\'evy-Fokker-Pl\-anck  equation}
\label{sec-fp}

The results presented in this section are coming from~\cite{gi}. We are looking for asymptotic  behaviour 
of the solution of a Fokker-Planck equation where the classical Laplacian is replaced with 
a L\'evy operator. Precisely, recalling that $\mathcal I$ is defined in~\eqref{operator},
we consider the L\'evy-Fokker-Planck equation 
\begin{equation}\label{eq:fp}
\partial_t u = \mathcal{I}[u] + \mathrm{div} (u \F) \quad x \in \R^d, t >0
\end{equation}
submitted to the initial condition
\begin{equation*} \label{eq:ci}
u(0,x) = u_0 (x) \quad x \in \R^d
\end{equation*}
where $u_0$ is nonnegative and in $L^1 (\R^d)$ and $F$ is a given proper force
for which there exists a nonnegative steady state (see below). 

\subsection{The $\Phi$-Entropy  and associated Fisher information}
\label{sec-ou-phi-fi}

In this subsection, we are interested in  the (time) derivative of the $\Phi$-entropy associated to the 
L\'evy-Fokker-Planck equation when a steady state is given. 
\begin{prop}
\label{prop:fisher}
Assume that there exists $u_\infty$, a steady state of~\eqref{eq:fp}, a positive solution of the equation :
\begin{equation} \label{eq:K}
\mathcal{I} [u_\infty] + \dive (u_\infty F) =0,
\end{equation}
such that $\int u_\infty dx=1$.  
  Assume that the initial condition $u_0$ is nonnegative and satisfies   
$\mathrm{Ent}^\Phi_{u_\infty}  \left( \frac{u_0}{u_\infty} \right)<\infty $.  

Then for any convex smooth function $\Phi:\R^+ \to \R$ and any $t\ge 0$,
the solution $u$ of \eqref{eq:fp} satisfies 
\begin{multline} 
\label{eq:decrease}
\forall t\geq 0,\quad
\frac{d}{dt}\mathrm{Ent}^\Phi_{u_\infty} ( v) = 
-\int\Phi''(v) \nabla v \cdot \sigma \nabla v\, u_\infty dx\\
- \int \int D_\Phi \left(v(t,x), v(t,x-z) \right) \nu (dz) u_\infty (x) dx
\end{multline}
where $v(t,x) =\frac{u(t,x)}{u_\infty (x)}$ and $\nu$ is the L\'evy measure 
appearing in the definition of the operator $\mathcal{I}$ and $D_\Phi$ is defined in~\eqref{eq-breg}.
\end{prop}
\medskip

In order to prove Proposition~\ref{prop:fisher}, since the $\Phi$-entropy involves the function 
$v(t,x)= \frac{u(t,x)}{u_\infty (x)}$,
its derivative makes appear $\partial_t v$ and it is natural to ask ourselves which partial differential 
equation $v$ satisfies. Using~\eqref{eq:K} one gets by a  simple computation 
\begin{eqnarray} \label{eq:ol}
\partial_t v &=& \frac1{u_\infty} \big( \mathcal{I} [ u_\infty v ] + \dive (u_\infty v \F) \big)
\nonumber\\
&=& \frac1{u_\infty} \big(\mathcal{I} [ u_\infty v ] - \mathcal{I} [ u_\infty  ]v \big) 
+ \F \cdot \nabla v=:Lv . 
\end{eqnarray}
In the case where $\mathcal{I} [u] = \Delta u$ 
(\textit{i.e.} $\sigma$ is the identity matrix, $b=0$ and $a=0$), 
Equation~\eqref{eq:ol} becomes
$$
\partial_t v = \Delta v - \F \cdot \nabla v
$$
and is known as the Ornstein-Uhlenbeck equation. This is the reason why we will refer to Equation~\eqref{eq:ol} 
as the L\'evy-Ornstein-Uhlenbeck equation.  We next give a simpler formulation for the L\'evy-Ornstein-Uhlenbeck operator. 
%
\begin{lem}[L\'evy-Ornstein-Uhlenbeck equation] \label{lem:ol}
If the integrodifferential operator on the right-hand side of \eqref{eq:ol} is denoted by $L$, 
we have for all smooth functions $w_1$ and $w_2$
$$
\int w_1 \; L w_2 \; u_\infty dx = \int (\check{\mathcal{I}} [w_1] 
- \F \cdot \nabla w_1) \; w_2 \; u_\infty dx
$$ 
 where  $\check{\mathcal{I}}$ is the L\'evy operator whose parameters are $(-b,\sigma,\check{\nu})$
 with $\check{\nu} (dz) = \nu (-dz)$.
This can be expressed by the formula: $L^* = \check{\mathcal{I}}  - \F \cdot \nabla $ 
where duality is understood with respect to the measure $u_\infty dx$.
\end{lem}
\begin{proof}
The main tool is the integration by parts for the operator $\mathcal I$, see equation~\eqref{eq-ipp}. 
For any smooth functions $u$, $v$ one gets 
$$
\int v\mathcal{I}[u]dx=\int u \check{\mathcal{I}} [v]dx.
$$ 
If $w_1$ and $w_2$ are two smooth functions on $\R^d$, then:
\begin{eqnarray*}
\int w_1 \; L w_2 \; u_\infty dx= \int w_1 \; ( \mathcal{I} [ u_\infty w_2 ] 
- \mathcal{I} [u_\infty] w_2 
+ u_\infty \F \cdot \nabla w_2) dx\\
 = \int w_2 \;  \check{\mathcal{I}} [ w_1 ] \; u_\infty dx - \int \mathcal{I} [u_\infty] w_1 w_2 dx
- \int \dive (u_\infty w_1 \F ) w_2 dx\\ 
= \int u_\infty (\check{\mathcal{I}} [w_1] 
- \F \cdot \nabla w_1)
w_2 dx. 
\end{eqnarray*}
\end{proof}

\begin{proof}[{Proof of Proposition~\ref{prop:fisher}.}]
By using Lemma~\ref{lem:ol} with $v=u/u_\infty$, we get:
\begin{eqnarray*}
\frac{d}{dt}\mathrm{Ent}^\Phi_{u_\infty} ( v)&=& \int \Phi'(v) \; \partial_t v \; u_\infty dx 
= \int \Phi'(v) \; Lv \; u_\infty dx\\
&=& \int \check{\mathcal{I}} [\Phi'(v)]\; v \; u_\infty dx -\int \F \cdot \nabla (\Phi'(v))  \; v \; u_\infty dx.
\end{eqnarray*}
If now one remarks that $r\Phi''(r) = (r \Phi' (r) -\Phi (r))'$, we get:
\begin{eqnarray*}
\frac{d}{dt}\mathrm{Ent}^\Phi_{u_\infty} ( v) &=&\int v \check{\mathcal{I}} [\Phi'(v)] u_\infty dx 
-\int \F \cdot \nabla (v\Phi'(v) - \Phi (v)) ) \; u_\infty dx \\
&=&\int v \check{\mathcal{I}} [\Phi'(v)] u_\infty  dx + \int 
\dive (u_\infty \F) (v\Phi'(v) - \Phi (v)) dx \\
&=&\int v \check{\mathcal{I}} [\Phi'(v)] u_\infty  dx 
- \int \mathcal{I} [u_\infty] (v\Phi'(v) - \Phi (v)) dx \\
&=&\int ( v \check{\mathcal{I}} [\Phi'(v)]  - \check{\mathcal{I}} [ v\Phi'(v) ] + \check{\mathcal{I}} [\Phi (v)] 
) u_\infty dx 
\end{eqnarray*}
\begin{eqnarray*}
\frac{d}{dt}\mathrm{Ent}^\Phi_{u_\infty} ( v)
&=&\int ( v {\check{\mathcal{I}}}_g [\Phi'(v)]  - {\check{\mathcal{I}}}_g [ v\Phi'(v) ] + \check{\mathcal{I}}_g [\Phi (v)] 
) u_\infty dx \\
& &+\int ( v \check{\mathcal{I}}_a [\Phi'(v)]  - \check{\mathcal{I}}_a [ v\Phi'(v) ] + \check{\mathcal{I}}_a [\Phi (v)] 
) u_\infty dx \\
 &=&-\int\Phi''(v) \nabla v \cdot \sigma\nabla v \;  u_\infty dx 
\\ & & \hspace{-1cm} +\int \int \bigg( v (x) (\Phi'(v(x+z)) - \Phi'(v(x)))
 - v(x+z) \Phi' (v (x+z)) \\
&&+ v (x) \Phi' (v(x))  + \Phi(v(x+z)) - \Phi (v(x))  ) \bigg) \check{\nu} (dz) \; u_\infty (x)dx
\end{eqnarray*}
\begin{eqnarray*}
\frac{d}{dt}\mathrm{Ent}^\Phi_{u_\infty} ( v) & =&-\int\Phi''(v) \nabla v \cdot \sigma\nabla v \; u_\infty dx\\
 && \quad -\int \int \bigg( \Phi (v(x)) -  \Phi(v(x+z))  -\Phi'(v(x+z)) \\
&& \qquad  \times  (v(x) -v (x+z)) \bigg)\,\check{\nu} (dz) \;  
u_\infty (x)dx.
\end{eqnarray*}
Then the definitions of the Bregman distance and of the L\'evy measure $\check{\nu}$ give
\begin{multline*}
\frac{d}{dt}\mathrm{Ent}^\Phi_{u_\infty} ( v)\! 
=\! -\int\Phi''(v) \nabla v \cdot\sigma \nabla v \; u_\infty dx\!\\
-\! \int \int D_\Phi (v(x), v(x-z))  
\nu (dz) \,u_\infty(x) dx.
\end{multline*}
\end{proof}

\subsection{Exponential decay of the L\'evy-Fokker-Planck equation}
\label{sec-convergence}
We give now assumptions such that there exists a steady state of the L\'evy-Fokker-planck equation. For this section we need  
to assume for that  the force is given by $F(x) = x$.

\begin{thm}[Exponential decay to equilibrium]
\label{thm:exp}
Assume that $F(x) = x$ and the operator $\mathcal I$ is the infinitesimal 
generator of a L\'evy process whose L\'evy measure is denoted by $\nu$. 
We assume that $\nu$ has a density  $N$  with respect to the Lebesgue measure and that $N$ satisfies
\begin{equation}
\label{eq:con1}
\int_{\R^d \setminus B} \ln |z| \; N(z) \; dz < +\infty
\end{equation}
where $B$ is the unit ball in $\R^d$.  

\begin{itemize}
\item
Then there exists a steady state $u_\infty$, \textit{i.e.} a nonnegative
solution of \eqref{eq:K} satisfying $\int u_\infty dx=1$. 
\item
If moreover $N$ is even and for all $z\in\R^d$,
\begin{equation}
\label{eq:con2}
\int_1^{+\infty} N  (s z) s^{d-1} ds \leq C N(z)
\end{equation}
for some constant $C\geq0$, 
then for any smooth convex function $\Phi$ such that condition~\eqref{cond:convex} is satisfied, 
the $\Phi$-entropy of the solution $u$ of \eqref{eq:fp}-\eqref{eq:ci} 
goes to 0 exponentially. Precisely, for any nonnegative initial 
datum $u_0$ such that 
$\mathrm{Ent}^\Phi_{u_\infty}  \left( \frac{u_0}{u_\infty} \right)<\infty,$
one gets:
\begin{equation}
\label{eq-ththm}
\forall t\ge0,\quad \mathrm{Ent}^\Phi_{u_\infty}  
\left( \frac{u(t)}{u_\infty} \right)\leq e^{-\frac{t}{C}}\mathrm{Ent}^\Phi_{u_\infty}  
\left( \frac{u_0}{u_\infty} \right)
\end{equation}
with $C$ appearing in \eqref{eq:con2}. 
\end{itemize}
\end{thm}

To prove the first part of the Theorem, the existence of the steady state, we need to state the following lemma. 
%
\begin{lem}\label{prop:steadystate}
Assume that  the L\'evy measure $\nu$ has a density $N$ with respect to the Lebesgue measure  
and that it satisfies~\eqref{eq:con1}.
There then exists a steady state $u_\infty$, \textit{i.e.} a 
solution of \eqref{eq:K}. Moreover, it is an
infinitely divisible measure whose characteristic exponent $A$ is defined by:
\begin{equation}\label{eq:defA}
A ( \xi ) =  - \xi\cdot  \sigma\xi + i b\cdot \xi 
+\int_0^1 a(s \xi) \frac{ds}s.
\end{equation}
Moreover, parameters of the characteristic exponent $A$ are  $(\sigma,b-b_A,N_\infty dx)$ where 
\begin{equation}\label{eq:defba}
b_A = \int \int_0^1 z \frac{(1-\tau^2) |z|^2}{(1+\tau^2 |z|^2)(1+|z|^2)}d\tau N (z) dz,
\end{equation}
and 
\begin{equation} \label{eq:NA}
N_\infty (z) = \int_1^\infty N (t z) t^{d-1} dt.
\end{equation}
Note that the L\'evy measure $\nu_\infty$ associated to the characteristic exponent $A$  has a density $N_\infty$
with respect to the Lebesgue measure.
\end{lem}
\begin{rem}
In the general case, Condition~\eqref{eq:con2} precisely says that $N_\infty \le C N$ which 
can be written  in terms of measures as follows: $\nu_\infty\leq C\nu$.
\end{rem}
\begin{proof}[Proof of Lemma~\ref{prop:steadystate}]
Let us start as in \cite{bk03} in section~3. At least formally, 
the Fourier transform $\hat{u}_\infty$ of any steady state $u_\infty$ satisfies 
$$
\psi(\xi)\hat{u}_\infty + \xi \cdot \nabla \hat{u}_\infty=0,
$$
where $\psi$ is the characteristic exponent of the L\'evy operator $\mathcal I$.  So that $\hat{u}_\infty = \mathrm{exp} (- A)$ with $A$ such that: 
$$
\nabla A (\xi)\cdot \xi = \psi(\xi).
$$
The solution of this equation is precisely given by \eqref{eq:defA}.
It is not clear that $A$ is well defined and is 
the characteristic exponent of an infinitely divisible measure; this is what
we prove next. This will imply  in particular that $\mathcal{F}^{-1} (\exp (-A))$ 
is a nonnegative function. 

Define the nonnegative $N_\infty $ 
by Equation~\eqref{eq:NA}. This integral of a nonnegative function is finite
since for any $R >0$, if $d\sigma$ denotes the uniform measure on the unit sphere
$S^{d-1}$ we get,
$$
 \int_{R}^\infty \int_{|D|=1} N (\tau D)\tau^{d-1}d\tau d\sigma (D)=  \int_{|y| \ge R} N (y) 
dy < +\infty. 
$$
We conclude that for any $r \ge R >0$ and almost every $D$ on the unit sphere (where the set of 
null measure depends only on $R$), 
$$
r^d N_\infty (rD) = \int_{r}^\infty  N (\tau  D) \tau^{d-1} d\tau < +\infty
$$
so that $N_\infty (z)$ is well-defined almost everywhere outside $B_R$. 
Choose now a sequence $R_n \to 0$ and conclude.

Let us define $I (r) =   \int_{|D|=1} N (rD) d\sigma (D)$ and $I_\infty=  \int_{|D|=1} N_\infty (rD) d\sigma (D)$ in an analogous way. 
The previous equality implies that:
$$
r^d I_\infty (r) = \int_r^{+\infty} I (\tau) \tau^{d-1} d\tau.
 $$
We conclude that:
\begin{eqnarray*}
\int_{|z| \le 1} |z|^2 N_\infty (z) dz &=& \int_0^1 I_\infty (r)  r^{d+1}dr =  \int_0^1 r \int_r^{+\infty} 
I (\tau) \tau^{d-1} d\tau dr \\
& = & \frac12 \int_{|z|\ge 1} N (z) dz + \frac12 \int_{|z| \le 1} |z|^2 N (z) dz  < +\infty\\
\int_{|z| \ge 1} N_\infty (z) dz &=& \int_1^{+\infty} I_\infty (r) r^{d-1} dr =  \int_1^{+\infty} \frac1r 
 \int_r^{+\infty}  I (\tau) \tau^{d-1} d\tau dr\\
& = &\int_{|z| \ge 1} \log |z| N (z) dz.  
\end{eqnarray*}
Hence we have  $\int \min (1,|z|^2) N_\infty (z) dz <+\infty$.
We conclude that it is a L\'evy measure. Now consider the associated characteristic exponent:
$$
\tilde{A} (\xi) =
i b\cdot \xi - \sigma \xi\cdot \xi
+ \int (e^{i z \cdot \xi} -1 - i (z \cdot \xi) \; h(z) ) N_\infty (z) dz.
$$
Now compute:
\begin{eqnarray*}
\tilde{A} (\xi) \!\!\! \!&-& i b\cdot \xi + \sigma \xi\cdot \xi \\
&= &\int \int_1^\infty 
(e^{i z \cdot  \xi} -1 - i (z \cdot \xi) \; h(z) ) N (sz) s^{d-1} ds \; dz \\
& =& \int_1^\infty \left\{ \int 
\left(e^{i \tilde{z} \cdot \frac{\xi}s} -1 - i \left(\tilde{z} \cdot  \frac\xi{s} \right) \; h\left( \frac{\tilde{z}}{s} \right) \right) 
N (\tilde{z})  d\tilde{z} \right \} \frac{ds}s \\
& =& \int_1^\infty a\left (\frac\xi{s} \right) \frac{ds}s -i \xi \cdot \int_1^\infty \left\{ \int
\frac{\tilde{z}}s  \left( h\left( \frac{\tilde{z}}{s} \right) - h (\tilde{z}) \right) 
N (\tilde{z})  d\tilde{z} \right \} \frac{ds}s \\
&=& A (\xi) -i \xi \cdot b_A 
\end{eqnarray*}
where $b_A$ is defined by \eqref{eq:defba}. Properties~\eqref{eq-def1} of the L\'evy measure $\nu$ imply 
that $b_A$ is well defined.  
We conclude that $A$ is the characteristic exponent
of an infinitely divisible law $u_\infty$ whose drift is $b-b_A$, whose 
Gaussian part is $\sigma$  and whose L\'evy measure is $N_\infty (z) dz$. 
\end{proof}

\begin{proof}[Proof of Theorem~\ref{thm:exp}]  The proof of the first part is exactly given by  Lemma~\ref{prop:steadystate}. 
 
We now turn to the second part of the theorem. 
Proposition~\ref{prop:fisher} gives for $t\geq0$,
\begin{eqnarray*}
\frac{d}{dt}\mathrm{Ent}^\Phi_{u_\infty}  
(v) = -\int\Phi''(v(t,\cdot)) \nabla v \cdot (t,\cdot) \cdot\sigma\nabla v (t,\cdot) u_\infty dx
\\ -\int \int D_\Phi ( v(t,x), 
v(t,x-z))\nu (dz)\, u_\infty dx \\
 = -\int\Phi''(v(t,\cdot)) \nabla v \cdot (t,\cdot)\cdot \sigma\nabla v (t,\cdot) u_\infty dx
\\ -\int \int D_\Phi ( v(t,x), 
v(t,x+z))\nu (dz)\, u_\infty dx
\end{eqnarray*}
where we used the fact that $\nu$ is even. 
It is now enough to prove the following inequality
\begin{eqnarray*}
\mathrm{Ent}^\Phi_{u_\infty}  \left( v \right) 
\leq C\int\Phi''(v(t,\cdot))\nabla v(t,\cdot)\cdot  \sigma\nabla v  u_\infty dx\\
+C \int \int D_\Phi \PAR{v(x),v(x+z)} \nu (dz)u_\infty(x) dx
\end{eqnarray*}
for some constant $C$ not depending on $v$ and Gronwall's lemma permits to conclude.
But this inequality is a direct consequence of \eqref{eq-ls} for the infinitely divisible
law $u_\infty$. 


\end{proof}
 \subsection{Examples}
We next discuss the condition we impose in order to get exponential decay, namely Condition~\eqref{eq:con2}. 
We point out that equality in this condition holds true only for $\alpha$-stable
operators and we give a necessary condition on the behaviour of the L\'evy measure 
at infinity if one knows that it decreases faster than $|x|^{-d}$.
\begin{prop} \label{prop:discuss}
\begin{itemize}
\item
Equality  $N_\infty=N/\lambda$ holds if and only if
$\psi$ is positively homogenous of index $\lambda\in(0,2]$, \textit{i.e.} 
$$\psi (t \xi) = t^{\lambda} \psi(\xi) \mbox{ for any } t>0, \xi \in \R^d.$$ 
In this case, we get $A=\psi/\lambda$ and $b_A=0$. Note that in the limit case 
$\lambda=2$, then  we  get $N_\infty=N/2=0$.
\item
If $|x|^d N (x) \to 0$ as $|x| \to +\infty$, then the densities $N$ and $N_\infty$ satisfy:
$$
N = - \mathrm{div} (x N_\infty).
$$
\item
In this case,  Condition~\eqref{eq:con2} is equivalent to:
\begin{equation}\label{eq:cn}
\left\{ \begin{array}{ll}
N_\infty (tx) \le N_\infty (x) t^{-d -1/C} \mbox{ if } t \ge 1 \\
N_\infty (tx) \ge N_\infty (x) t^{-d -1/C} \mbox{ if } 0<t \le 1 
\end{array} \right.
\end{equation}
\end{itemize}
\end{prop}
\begin{proof}
The first item simply follows from the definition of $A$. 

Let us first prove the second item. 
\begin{eqnarray*}
N (x) &=& - (\lim_{t\to +\infty}  t^d N (tx))+
1^d N (1 \times x) = - \int_1^{+\infty} \frac{d}{dt} (t^d N (tx)) dt \\
&=& -d N_\infty (x) - x\cdot \nabla N_\infty (x) = - \mathrm{div} (x N_\infty).
\end{eqnarray*}
To prove the third item, use the first one to rewrite \eqref{eq:con2} as follows:
$$
x \cdot\nabla N_\infty (x) + (d+1/C) N_\infty (x) \le 0.
$$
Integrate over $[1,t]$ for $t\geq 1$ and $[t,1]$ for $t\leq 1$ to get the result. 
\end{proof}
\begin{ex}
In $\R$, the L\'evy measure $\frac1{|z|}e^{-|z|}$ does not satisfy Condition~\eqref{eq:con2}.  
Indeed, it is equivalent to:
$$
\int_1^{+\infty} \frac{ e^{ - |x| (s-1)}}{s} ds \le C
$$
and the monotone convergence theorem implies that the left hand side of this inequality
goes to $+\infty$ as $|x| \to +\infty$. 
\end{ex}

\begin{eremer}
The first author sincerely thanks H. Ouerdiane for his hospitality during the autumn in Hammamet. 
\end{eremer}


\end{document}